\documentclass[11pt]{amsart}

\usepackage{amsmath,amsthm,amscd,euscript}
\newtheorem{theorem}{Theorem}

\theoremstyle{remark}

\theoremstyle{definition}

\newtheorem{problem}{Problem}
\newtheorem{conjecture}[problem]{Conjecture}

\title{Approximating reals  by rationals of the form $a/b^2$.}
\author{Oleg~Karpenkov}
\date{30 October 2006}
\thanks{Partially supported
by NWO-RFBR 047.011.2004.026 (RFBR 05-02-89000-NWO\_a) grant, by
RFBR SS-1972.2003.1 grant, by RFBR 05-01-02805-CNRSL\_a grant,
and by RFBR grant 05-01-01012a.}


\email[Oleg Karpenkov]{karpenk@mccme.ru}

\address{Mathematisch Instituut, Universiteit Leiden,
P.O. Box 9512, 2300 RA Leiden, The Netherlands}

\begin{document}
\input epsf

\begin{abstract}
In this note we formulate some questions in the study of
approximations of reals by rationals of the form $a/b^2$ arising
in theory of Shr\"odinger equations. We hope to attract attention
of specialists to this natural subject of number theory.
\end{abstract}

\maketitle

\sloppy \normalsize

\section{Introduction}

{\it Some background.} In this note we formulate some questions
in the study of approximations of reals by rationals of the form
$a/b^2$ arising in theory of Shr\"odinger equations
(see~\cite{Cra} and~\cite{Kri} for further information). We hope
to attract attention of specialists to this natural subject of
number theory. Good references to theory of approximations by
arbitrary algebraic numbers are for instance~\cite{Wir},
\cite{Spr} and~\cite{Tis}, especially for approximations by
quadratic irrationals, see~\cite{Dav}. A metric approach to the
study (in a more general situation) was proposed by~\cite{Sch}
and further developed by M.~Weber in~\cite{Web}, V.~Beresnevich,
M.~Dodson, S.~Kristensen, and J.~Levesley in~\cite{Ber} and other
works. This approach is a good test of the proposed problems,
nevertheless it does not give the answers. Some upper bound
estimates were made by A.~Zaharescu in~\cite{Zah}.

\section{Questions}

\vspace{1mm}

We start with formulation of one of the main results in classical
approximation theory of reals by rationals (see~\cite{Khi} for the
proofs).

\vspace{1mm}

\begin{theorem} {\bf I.} For any reals $\alpha$ and $c\ge 1/\sqrt{5}$ there
exists an infinite number of integer solutions $(a,b)$, $b>0$ for
the following inequality
$$
\left|\alpha - \frac{a}{b}\right|< \frac{c}{b^2}.
$$

{\bf II.} Let  $\alpha$ be the Golden Ratio $($i.~e.
$\alpha=(\sqrt{5}+1)/2$$)$. Then for any $c< 1/\sqrt{5}$ the
inequality of item~I. has only finitely many solutions. \qed
\end{theorem}

\vspace{1mm}

Similar results for the approximations by rationals of the form
$a/b^2$ are not known. One of the reasons of that is the
following: lattice geometry of continued fractions corresponding
to the approximations by rationals can not be naturally expanded
to the case of rationals $a/b^2$.

\begin{problem}
Find a good generalization for geometry of numbers to the case of
$a/b^2$-approximations.
\end{problem}

\vspace{1mm}

Let us give some known estimates for the case of
$a/b^2$-approximations. The lower estimate seems to be quite
precise.

\begin{theorem}\label{t1}
For any positive $\varepsilon$ and for any real $\alpha$ there
exist a positive constant $c=c(\varepsilon)$, such that the
following inequality does not have integer solutions:
$$
\left|\alpha
-\frac{a}{b^2}\right|<\frac{c}{b^3\ln^{1+\varepsilon}(b)}.
$$
\end{theorem}

The proof of a more general statement is given by I.~Borosh and
A.~S.~Fraenkel in~\cite{Bor}. We suppose that the logarithm in
the formula can be eliminated.

All known proofs of previous theorem are general, and do not give
the examplpes of badly approximable reals, like it was
$(\sqrt{5}+1)/2$ for the case of approximations by $a/b$. So the
following problem is actual here.

\begin{problem}
Find any particular example of $\alpha$ that satisfies the
condition of Theorem~\ref{t1}.
\end{problem}

The following estimate for the upper bound case is known.

\begin{theorem}\label{t2}{\bf A.~Zaharescu~\cite{Zah}.}
For any real $\alpha$ and any positive $\theta<2/3$ there exists
infinitely many solutions of the following inequality:
$$
\left|\alpha -\frac{a}{b^2}\right|<\frac{1}{b^{2+\theta}\xi(b)}.
$$
\end{theorem}

As one can see there is a gap between upper and lover estimates
for the baddly approximable reals by $a/b^2$-rationals. The
results~\cite{Bor} for almost all reals and numerical experiments
support the following classical conjecture.

\begin{conjecture}
For any real $\alpha$ there exists a constant $c(\alpha)$ such
that the inequality:
$$
\left|\alpha -\frac{a}{b^2}\right|<\frac{c(\alpha)}{b^3}
$$
has infinitely many integer solutions $(a,b)$.
\end{conjecture}

We conclude this note with the following problem wich is supposed
to be close to the subject.

\begin{problem}
Find the estimates for the upper and lower bounds of the ``best''
approximations of reals by rationals of the form $z/p$, where $z$
is integer, and $p$ is prime or unity.
\end{problem}

{\bf Acknowledgement.} The author is grateful to E.~S\'er\'e,
W.~Craig, H.~W.~Lenstra, J.-H.~Evertse, S.~Kristensen, and Nigel
Watt, for help with collecting the information and comments, and
Mathematisch Instituut of Universiteit Leiden for the hospitality
and excellent working conditions.


\begin{thebibliography}{99}

\bibitem{Ber} V.~Beresnevich, M.~Dodson, S.~Kristensen and J.~Levesley,
{\it An inhomogeneous wave equation and non-linear Diophantine
approximation}, preprint, March 2006, 18~p.

\bibitem{Bor} I.~Borosh, A.~S.~Fraenkel, {\it A generalization of JarnМk's theorem
on Diophantine approximations}, Indag. Math. 34 (1972),
pp.~193--201.

\bibitem{Cra}
W.~Craig, {\it Probl\`emes de petits diviseurs dans les
\'equations aux d\'eriev\'ees partielles}, Panoramas et
synth\`eses, 9, SMF, Paris 2000.


\bibitem{Dav} H.~Davenport, W.~M.~Schmidt, {\it Approximation to real numbers
by quadratic irrationals}, Acta Arith., 13(1967/1968),
pp.~169--176.

\bibitem{Kri} S.~Kristensen, {\it Diophantine approximation and the solubility of
the Schrodinger equation}, Phys. Lett. A, 314(1-2), 2003,
pp.~15-18.

\bibitem{Khi}
A.~Ya.~Khinchin, {\it Continued fractions}, Moscow, FIZMATGIZ,
3.~ed., (1961); English translation University of Chicago Press,
1961.

\bibitem{Sch} W.~M.~Schmidt, {\it Metrical theorems on fractional parts of
sequences}, Trans. Amer. Math. Soc., 110(1964), pp.~493--518.

\bibitem{Spr}
V.~G.~Sprind\v zuk, {\it A proof of Mahler's conjecture on the
measure of the set of $S$-numbers} (Russian), Izv. Akad. Nauk SSSR
Ser. Mat. 29(1965), pp.~379--436.

\bibitem{Tis}
K.~I.~Tishchenko, {\it On approximation to real numbers by
algebraic numbers}, Acta Arith., 94(2000), no.~1, pp.~1--24.

\bibitem{Web}
M.~Weber, {\it Some examples of application of the metric entropy
method}, Acta Math. Hungar. 105(2004), no.~1-2, pp.~39--83.

\bibitem{Wir}
E.~Wirsing, {\it Approximation mit algebraischen Zahlen
beschrankten Grades} (German), J. Reine Angew. Math., 206(1960),
pp.~67--77.

\bibitem{Zah}, A.~Zaharescu, {\it Small values of $n\sp 2\alpha\pmod
1$}, Invent. Math., 121(1995), no.~2, pp.~379--388.

\end{thebibliography}
\end{document}